\documentclass{amsart}

\newtheorem*{theorem}{Theorem}

\begin{document}

\title[A proof of the Steiner-Lehmus Theorem in hyperbolic geometry]{A trigonometric proof of the Steiner-Lehmus Theorem in hyperbolic geometry}

\author{Keiji Kiyota}
\address{Graduate School of Integrated Basic Sciences, Nihon University,
3-25-40 Sakurajosui, Setagaya-ku, Tokyo 156-8550, Japan}
\email{ta15039@educ.chs.nihon-u.ac.jp}

\begin{abstract}
We give a trigonometric proof of the Steiner-Lehmus Theorem in hyperbolic geometry. Precisely we show that if two internal bisectors of a triangle on the hyperbolic plane are equal, then the triangle is isosceles.
\end{abstract}



\date{\today}

\maketitle

\section{Introduction}
In 1844 \cite{JS}, Steiner gave the first proof of the following theorem. If two internal bisectors of a triangle on the Euclidean plane are equal, then the triangle is isosceles. This had been originally asked by Lehmus in 1840, and now is called the Steiner-Lehmus Theorem. Since then, wide variety of proofs have been given by many people over 170 years. At present, at least 80 different proofs exist. See \cite{PSS}. For example, in 2008, Hajja gave a short trigonometric proof in \cite{MH}. On the other hand, several proofs of this theorem in hyperbolic geometry ware given in \cite{NS}, \cite{CB} and \cite{OD}. In this paper, we give a simple trigonometric proof in hyperbolic geometry based on the way of Hajja. 

 \section{Steiner-Lehmus Theorem}
\begin{theorem}
  If two internal bisectors of a triangle on the Hyperbolic plane are equal, then the triangle is isosceles.
\end{theorem}

\begin{figure}[htb]
{\unitlength 0.1in%
\begin{picture}(21.3000,17.3000)(6.7000,-20.0000)%
%
\special{pn 8}%
\special{pa 1800 400}%
\special{pa 800 2000}%
\special{fp}%
\special{pa 800 2000}%
\special{pa 2800 2000}%
\special{fp}%
\special{pa 2800 2000}%
\special{pa 1800 400}%
\special{fp}%
%
\special{pn 8}%
\special{pa 800 2000}%
\special{pa 2310 1200}%
\special{fp}%
\special{pa 1310 1200}%
\special{pa 2790 2000}%
\special{fp}%
\put(17.4000,-4.0000){\makebox(0,0)[lb]{$A$}}%
\put(6.7000,-21.1000){\makebox(0,0)[lb]{$B$}}%
\put(27.7000,-21.1000){\makebox(0,0)[lb]{$C$}}%
\put(23.1000,-12.1000){\makebox(0,0)[lb]{$B'$}}%
\put(11.3000,-12.1000){\makebox(0,0)[lb]{$C'$}}%
\put(10.3000,-20.5000){\makebox(0,0)[lb]{$\beta$}}%
\put(25.2000,-19.9000){\makebox(0,0)[lb]{$\gamma$}}%
\put(21.0000,-8.1000){\makebox(0,0)[lb]{$u$}}%
\put(26.1000,-16.1000){\makebox(0,0)[lb]{$U$}}%
\put(9.3000,-16.1000){\makebox(0,0)[lb]{$V$}}%
\put(14.0000,-8.2000){\makebox(0,0)[lb]{$v$}}%
\end{picture}}%
\label{1}
 \caption{}
\end{figure}

 \begin{proof}
 We consider a triangle $ABC$ on the hyperbolic plane. See Figure 1. Let $B'$ be intersection of the side $AC$ and the internal bisectors of the angle $B$. Let $C'$ be the intersection of the side $AB$ and the internal bisector of the angle $C$. Then $BB'$ and $CC'$ are the internal bisectors of the angles $B$ and $C$. Let $a , b$ and $c$ be the lenghs of the opposite sides of the angles $A , B$ and $C$ respectively. We set $\beta = B/2$, $\gamma = C/2$, $u = AB'$, $U = B'C$, $v = AC'$, and $V = C'B$ .\par
 We apply the sines theorem in hyperbolic geometry to the triangles $ABC$, $BCC'$, $ACC'$, $CBB'$ and $ABB'$ respectively, then we have the following.
  \begin{align}
  \frac{\sinh{a}}{\sin{A}} = \frac{\sinh{b}}{\sin{2 \beta}} = \frac{\sinh{c}}{\sin{2 \gamma}}\\
  \frac{\sinh{CC'}}{\sin{2 \beta}} = \frac{\sinh{V}}{\sin{ \gamma}}\\
  \frac{\sinh{CC'}}{\sin{A}} = \frac{\sinh{v}}{\sin{ \gamma}}\\
  \frac{\sinh{BB'}}{\sin{2 \gamma}} = \frac{\sinh{U}}{\sin{ \beta}}\\
  \frac{\sinh{BB'}}{\sin{A}} = \frac{\sinh{u}}{\sin{ \beta}}
  \end{align}
  We assume $BB' = CC'$ and $C > B$, and lead to contradiction. Since the sum of the interior angles in a hyperbolic triangle is less than $\pi$, we have $B < C < \frac{\pi}{2}$, and so, $\sin{B} < \sin{C}$.
  In the following, we evaluate the magnitude relationship of $u , v$ and $U , V$ respectively.\\
 By (2) and (4), 
 \begin{align*}
  \frac{\sinh{V}}{\sin{ \gamma}}\sin{2 \beta} &= \frac{\sinh{U}}{\sin{ \beta}}\sin{2 \gamma} \\
  \frac{\sinh{U}}{\sinh{V}} &= \frac{\sin{\beta}}{\sin{\gamma}}\frac{\sin{2 \beta}}{\sin{2 \gamma}} 
 \end{align*}
  By (1), we get $\frac{\sin{2 \beta}}{\sin{2 \gamma}} = \frac{\sinh{b}}{\sinh{c}}$, so we have the following.
 \[
  \frac{\sinh{U}}{\sinh{V}}  = \frac{\sin{\beta}}{\sin{\gamma}}\frac{\sinh{b}}{\sinh{c}}
 \]
  Becaue of $\frac{\sin{\beta}}{\sin{\gamma}} < 1$ and $\frac{\sinh{b}}{\sinh{c}} < 1$,  we get $\frac{\sin{\beta}}{\sin{\gamma}}\frac{\sinh{b}}{\sinh{c}} < 1$.Then we have $\sinh{U} < \sinh{V}$. Since the hyperbolic sine function is monotonically increasing, we conclude $U < V$.\\
Similarly, by (3) and (5),
  \begin{align*}
  \frac{\sinh{v}}{\sin{ \gamma}} &= \frac{\sinh{u}}{\sin{ \beta}} \\
  \frac{\sinh{u}}{\sinh{v}} &= \frac{\sin{ \beta}}{\sin{ \gamma}}  < 1
  \end{align*}
  Therefore we get $\sinh{u} < \sinh{v}$, that is, $u < v$.\par
 Now let us consider the ratio and difference of $\frac{\sinh{b}}{\sinh{u}}$ and $\frac{\sinh{c}}{\sinh{v}}$. First we consider the ratio.
  \[
 \frac{\sinh{b}}{\sinh{u}} \Bigl/ \frac{\sinh{c}}{\sinh{v}} = \frac{\sinh{b}}{\sinh{u}}\frac{\sinh{v}}{\sinh{c}} = \frac{\sinh{b}}{\sinh{c}}\frac{\sinh{v}}{\sinh{u}} 
 \]
 We have the following by (1).
 \[
\frac{\sinh{b}}{\sinh{c}}\frac{\sinh{v}}{\sinh{u}} = \frac{\sin{2 \beta}}{\sin{2 \gamma}}\frac{\sin{ \gamma}}{\sin{ \beta}}
 \]
 Here, we apply the double-angle formula to $\sin{2 \beta}, \sin{2 \gamma}$ respectively.
\[
\frac{\sin{2 \beta}}{\sin{2 \gamma}}\frac{\sin{ \gamma}}{\sin{ \beta}}= \frac{2\sin{ \beta}\cos{ \beta}}{2\sin{ \gamma}\cos{ \gamma}} \frac{\sin{ \gamma}}{\sin{ \beta}}= \frac{\cos{ \beta}}{\cos{ \gamma}}
\]
By assumption $\beta < \gamma$, we have $\cos{\beta} > \cos{\gamma}$. So $\frac{\cos{ \beta}}{\cos{ \gamma}} > 1$.
Therefore we get the following result.
\begin{equation}
 \frac{\sinh{b}}{\sinh{u}} > \frac{\sinh{c}}{\sinh{v}}
\end{equation}
Next we consider the difference.
 \[
    \frac{\sinh{b}}{\sinh{u}}-\frac{\sinh{c}}{\sinh{v}} = \frac{\sinh{(U + u)}}{\sinh{u}} - \frac{\sinh{(V + v)}}{\sinh{v}}
 \]
 We apply the sum formula to $\sinh{(U + u)}$ and $\sinh{(V +v)}$ respectively.
 \begin{align*}
 \frac{\sinh{(U + u)}}{\sinh{u}} - \frac{\sinh{(V + v)}}{\sinh{v}} &= \frac{\sinh{U}\cosh{u} + \cosh{U}\sinh{u}}{\sinh{u}} - \frac{\sinh{V}\cosh{v} + \cosh{V}\sinh{v}}{\sinh{v}} \\
									   &= \frac{\sinh{U}}{\sinh{u}}\cosh{u} + \cosh{U} - \frac{\sinh{V}}{\sinh{v}}\cosh{v} + \cosh{V}
 \end{align*}
 By (4), (5) and (2), (3), $\frac{\sinh{U}}{\sinh{u}} = \frac{\sin{A}}{\sin{2 \gamma}}$ and $\frac{\sinh{V}}{\sinh{v}} = \frac{\sin{A}}{\sin{2 \beta}}$ hold, and so, we have the following.
\[
\frac{\sinh{U}}{\sinh{u}}\cosh{u} + \cosh{U} - \frac{\sinh{V}}{\sinh{v}}\cosh{v} + \cosh{V} = \frac{\sin{A}}{\sin{2 \gamma}}\cosh{u} + \cosh{U} - \frac{\sin{A}}{\sin{2 \beta}}\cosh{v} + \cosh{V}
\]
Moreover we get the following by (1). 
\[
\frac{\sin{A}}{\sin{2 \gamma}}\cosh{u} + \cosh{U} - \frac{\sin{A}}{\sin{2 \beta}}\cosh{v} + \cosh{V} = \frac{\sinh{a}}{\sinh{c}}\cosh{u} + \cosh{U} - \frac{\sinh{a}}{\sinh{b}}\cosh{v} - \cosh{V}
\]
By $\sinh{c} > \sinh{b}$, we have $\frac{\sinh{a}}{\sinh{b}} > \frac{\sinh{a}}{\sinh{c}}$. And $\cosh{v} > \cosh{u}$ and $\cosh{V} > \cosh{U}$ by $u < v$ and $U < V$. Therefore we get the following.
\[
\frac{\sinh{a}}{\sinh{c}}\cosh{u} + \cosh{U} - \frac{\sinh{a}}{\sinh{b}}\cosh{v} - \cosh{V} < 0
\]
Eventually we conclude the following result.
\begin{equation}
 \frac{\sinh{b}}{\sinh{u}} < \frac{\sinh{c}}{\sinh{v}}
\end{equation} 
 A contradiction is led by (6) and (7).
 \end{proof}
 

\end{document}